# M/M/1 Queueing System with Non-preemptive Priority

Zhao Guo-xi

Department of Mathematics, Xinxiang University, Xinxiang 453000,China

Hu Qi-Zhou

Faculty of Science, Hefei University of Technology, Hefei 230009,China

**Abstract:** In this paper, the performance of non-preemptive M/M/1 queueing system with two priority is analyzed. By using complementary variable method to make vector Markov process and analyzing the state-change equations of the queueing system, the generating function of two kinds of customers' length   distribution are derived under non-preemptive priority .Through further discussion, the probability of the server that it is working or free and average length of two kinds of customers are also derived.

## 1.  Introduction

In order to offer different quality of service for different kinds of customers, we often control a queueing system by priority mechanism. This phenomenon is common in practice. For example, in telecommunication transfer protocol, for guaranteeing different layers service for different customers, priority classes control may appears in header of IP package, or in ATM cell. Priority control is also wildly used in production practice, transportation management, etc.

Familiar priority control involves preemptive priority and non-preemptive priority. Considering a queueing system with two kinds of customers, when a customer of first class arrives the server, he finds the server is serving for a customer of second class, he squeezes the customer -in-service out and receives service at once. At the same time, customers belonging to same class obey the FCFS discipline, this mechanism is called preemptive priority queueing. If the arriving customer of first class finds the server is serving for a customer of second class should he wait a period until the customer -in-service finishes it's service, then begins receiving service, customers belonging to same kind obey the FCFS discipline, this mechanism is called non-preemptive priority queueing. Obviously, first class of customers have high priority over the second class of customers.

Cohen. J. W analyzed a two-class M/G/1 queueing system with non-preemptive priority [3] (454-460), but conclusion in [3] is just an approximate result. Miller investigated a M/M/1 queueing system with non-preemptive priority [6], Edward P. C. Kao consider a M/M/N queueing system of two types of customers under preemptive priority [9], matrix-geometry method is a efficient method to solve priority problem, but the process of computing is very complicated. [4] and [5] all research the non-preemptive priority queueing system, finding fitting time point and formatting imbedded markov chains is a common method, though the non-preemptive priority problem is too complicated, it's hard to get precise solution finally. Bong Dae Choi [8] investigated M/M/1 queueing system with deadline until the end of service, looking the state of server as a variable and formatting a vector markov process [2], the question may become simple. In this article, we use complementary variable method too, through formatting vector markov process [2] and analyzing the state-change equations of the M/M/1 queueing system with non-preemptive priority regulation, we derive the generating function of two kinds of customers' length distribution and average length of two kinds of customers.

## 2. System description
### 2.1 Model desription

We consider a single server queueing system serving two types of customers: class-1 (or first class) and class-2 (or second class) customers. The arrival process for both classes is state independent, we base our analysis on the following assumption:
1. Arrivals for both classes are independent stationary poission processes, importation parameters are $\lambda_1, \lambda_2$ respectively.

2. Service time is same expontially distribution, average service time is $1/\mu$.

3. Class-1 customers have non-preemptive priority over class-2 customers, we assume buffer capacity is infinite.

### 2.2 State-space description

$N_i(t)$ $(i=1,2)$ means the number of class-$i$ customer present in the system at time $t$.

Obviously customers' length process $\{N_1(t), N_2(t)\}$ has not the property of no-memory, so it isn't markov process. Varible $\xi = 0,1,2$ respectively represents the three states that the server is free or it is serving class-1 customers or class-2 customers is in service. Then vector process $N(t) = \{(N_1(t), N_2(t), \xi) : t \geq 0\}$ is markov process, state-space can be wrote as following structure: （000）;（101）（111）（121）（131）…，（201）(211)(221)…,…; (012)(022)(032)…, (112)(122)(132)…,…

## 3. Stationary system description
### 3.1 Stationary state-change equation

If the two classes customers are saw as the same kind, this system becomes an ordinary M/M/1queueing system .At the same time, input process is a synthesized poisson process with parameter $\lambda_1 + \lambda_2$, stationary condition is average service strength $\rho = (\lambda_1 + \lambda_2)/\mu < 1$.

When $\rho < 1$, Denoting $p_{ijk} = \lim_{t \to \infty} \Pr\{N_1(t) = i, N_2(t) = j, \xi = k\}$, stationary state-change equations can be wrote as follows:

$$p_{000}(\lambda_1 + \lambda_2) = p_{101}\mu + p_{012}\mu \qquad (1´)$$

$$p_{012}(\lambda_1 + \lambda_2 + \mu) = p_{111}\mu + p_{022}\mu + p_{000}\lambda_2 \tag{2´}$$

$$p_{0j2}(\lambda_1 + \lambda_2 + \mu) = p_{1j1}\mu + p_{0(j+1)2}\mu + p_{0(j-1)2}\lambda_2 \tag{3´}$$

$$p_{112}(\lambda_1 + \lambda_2 + \mu) = p_{012}\lambda_1 \tag{4´}$$

$$p_{1j2}(\lambda_1 + \lambda_2 + \mu) = p_{0j2}\lambda_1 + p_{1(j-1)2}\lambda_2 \tag{5´}$$

$$p_{i12}(\lambda_1 + \lambda_2 + \mu) = p_{(i-1)12}\lambda_1 \tag{6´}$$

$$p_{ij2}(\lambda_1 + \lambda_2 + \mu) = p_{(i-1)j2}\lambda_1 + p_{i(j-1)2}\lambda_2 \tag{7´}$$

$$p_{101}(\lambda_1 + \lambda_2 + \mu) = p_{201}\mu + p_{112}\mu + p_{000}\lambda_1 \tag{8´}$$

$$p_{111}(\lambda_1 + \lambda_2 + \mu) = p_{211}\mu + p_{122}\mu + p_{101}\lambda_2 \tag{9´}$$

$$p_{1j1}(\lambda_1 + \lambda_2 + \mu) = p_{2j1}\mu + p_{1(j+1)2}\mu + p_{1(j-1)1}\lambda_2 \tag{10´}$$

$$p_{i01}(\lambda_1 + \lambda_2 + \mu) = p_{(i+1)01}\mu + p_{i12}\mu + p_{(i-1)01}\lambda_1 \tag{11´}$$

$$p_{i11}(\lambda_1 + \lambda_2 + \mu) = p_{(i+1)11}\mu + p_{i22}\mu + p_{(i-1)11}\lambda_1 + p_{io1}\lambda_2 \tag{12´}$$

$$p_{ij1}(\lambda_1 + \lambda_2 + \mu) = p_{(i+1)j1}\mu + p_{i(j+1)1}\mu + p_{(i-1)j1}\lambda_1 + p_{i(j-1)1}\lambda_2 \tag{13´}$$

### 3.2 stationary-equations solving

In order to solve stationary-equations above, we define some marks as follows:

$$F_i^1(z_2) = \sum_{j=0}^{\infty} p_{ij1} z_2^j, i \geq 1 \quad ; \qquad F_i^2(z_2) = \sum_{j=1}^{\infty} p_{ij2} z_2^j, i \geq 0;$$

$$F^1(z_1, z_2) = \sum_{i=1}^{\infty} F_i^1(z_2) z_1^i ; \qquad F^2(z_1, z_2) = \sum_{i=0}^{\infty} F_i^2(z_2) z_1^i, i \geq 0$$

Obviously, united-distribution generating function of the two classes of customers' length can be written as:

$$F(z_1, z_2) = F^1(z_1, z_2) + F^2(z_1, z_2) + p_{000}$$

From expression (2´)、(3´) and (1´), we have:

$$[\lambda_1 + \lambda_2(1 - z_2) + \mu(1 - \frac{1}{z_2})]F_0^2(z_2) = \mu F_1^1(z_2) - p_{000}(\lambda_1 + \lambda_2 - \lambda_2 z_2) \tag{1}$$

From expression (4´) and (5´), we derive:

$$(\lambda_1 + \lambda_2(1 - z_2) + \mu)F_1^2(z_2) = \lambda_1 F_0^2(z_2) \tag{2}$$

From expression （6´） and （7´），we get:

$$[\lambda_1 + \lambda_2(1-z_2) + \mu]F_i^2(z_2) = \lambda_1 F_{i-1}^2(z_2) \quad (3)$$

From expression (8´)、(9´) and （10´），we have：

$$[\lambda_1 + \lambda_2(1-z_2) + \mu]F_1^1(z_2) = \mu F_2^1(z_2) + \frac{\mu}{z_2}F_1^2(z_2) + p_{000}\lambda_1 \quad (4)$$

From expression （11´）、（12´） and （13´），we have：

$$[\lambda_1 + \lambda_2(1-z_2) + \mu]F_i^1 = \mu F_{i+1}^1(z_2) + \frac{\mu}{z_2}F_i^2(z_2) + \lambda_1 F_{i-1}^1(z_2) \quad (5)$$

Deducing from expression (1)、(2)and (3) results in：

$$[\lambda_1(1-z_1) + \lambda_2(1-z_2) + \mu]F^2(z_1,z_2) = [\lambda_1 + \lambda_2(1-z_2) + \mu]F_0^2(z_2)$$

That gives:

$$F^2(z_1,z_2) = [\lambda_1 + \lambda_2(1-z_2) + \mu]F_0^2(z_2) / [\lambda_1(1-z_1) + \lambda_2(1-z_2) + \mu] \quad (6)$$

From expression (4) and (5)：

$$[\lambda_1(1-z_1) + \lambda_2(1-z_2) + \mu(1-\frac{1}{z_1})]F^1(z_1,z_2) =$$

$$\frac{\mu}{z_2}F^2(z_1,z_2) - [\mu F_1^1(z_2) + \frac{\mu}{z_2}F_0^2(z_2)] + \lambda_1 z_1 p_{000}$$

Substitute expression (1) and (6) into upper expression results in:

$$[\lambda_1(1-z_1) + \lambda_2(1-z_2) + \mu(1-\frac{1}{z_1})]F^1(z_1,z_2) =$$

$$[\lambda_1 + \lambda_2(1-z_2) + \mu](\frac{\mu}{z_2}\frac{1}{\lambda_1(1-z_1) + \lambda_2(1-z_2) + \mu} - 1)F_0^2(z_2) - [\lambda_1(1-z_1) + \lambda_2(1-z_2)]p_{000}$$

Thus

$$F^1(z_1,z_2) = \frac{[\lambda_1 + \lambda_2(1-z_2) + \mu](\frac{\mu}{z_2}\cdot\frac{1}{\lambda_1(1-z_1) + \lambda_2(1-z_2) + \mu} - 1)F_0^2(z_2) - [\lambda_1(1-z_1) + \lambda_2(1-z_2)]p_{000}}{\lambda_1(1-z_1) + \lambda_2(1-z_2) + \mu(1-\frac{1}{z_1})} \quad (7)$$

Considering the following equation (the left side is denominator of expression (7) ):

$$\lambda_1(1-z_1) + \lambda_2(1-z_2) + \mu - \frac{\mu}{z_1} = 0$$

That is：$\lambda_1 z_1^2 - [\lambda_1 + \lambda_2(1-z_2) + \mu]z_1 + \mu = 0 \quad (8)$

When $|z_1| = 1, |z_2| < 1$， we have the following inequality:

$$|\lambda_1 z_1^2 + \mu| \leq \lambda_1 + \mu < |\lambda_1 + \lambda_2(1-z_2) + \mu| = |(\lambda_1 + \lambda_2(1-z_2) + \mu)z_1|$$

From Rouche theorem, expression (8) has unique root in the area of $|z_2| < 1$, we may as well denote the root of equation (8) as $f(z_2)$. For $F^1(z_1, z_2)$ converges on the condition of $|z_1| < 1$ and $|z_2| < 1$, obviously $f(z_2)$ is the root of numerator of (7), then the expression as follow is correct:

$$[\lambda_1 + \lambda_2(1-z_2) + \mu](\frac{\mu}{z_2} \cdot \frac{1}{\lambda_1(1-f(z_2)) + \lambda_2(1-z_2) + \mu} - 1)F_0^2(z_2) - [\lambda_1(1-f(z_2)) + \lambda_2(1-z_2)]p_{000} = 0$$

From which:

$$F_0^2(z_2) = \frac{[\lambda_1(1-f(z_2)) + \lambda_2(1-z_2) + \mu] \cdot [\lambda_1(1-f(z_2)) + \lambda_2(1-z_2)]}{[\frac{\mu}{z_2} - (\lambda_1(1-f(z_2)) + \lambda_2(1-z_2) + \mu)] \cdot [\lambda_1 + \lambda_2(1-z_2) + \mu]} p_{000} \quad (9)$$

Then (6) and (7) can be expressed formulas including $p_{000}$

Substitute $z_1 = f(z_2)$ into expression (8), We have:

$$\lambda_1 f(z_2)^2 - [\lambda_1 + \lambda_2(1-z_2) + \mu]f(z_2) + \mu = 0 \quad (10)$$

When $z_2$ converges to 1, expression (10) can turn into the following expression:

$$\lambda_1 f(1)^2 - (\lambda_1 + \mu)f(1) + \mu = 0 \quad (11)$$

After solving equation (11), we easily get $f(1) = 1$ or $f(1) = \frac{\mu}{\lambda_1}$ (the latter must be deleted because $\mu/\lambda > 1$). As $f(z_2)$ is the root of square equation concerning $z_1$, $f(z_2)$ must be elementary function, it's differentiable. When $z_2$ converges to 1, derivating both sides of expression (10), we have $f'(1) = \frac{\lambda_2}{\mu - \lambda_1}$.

Substitute $f'(1) = \frac{\lambda_2}{\mu - \lambda_1}$ into expression (9), when $z_2$ converges to 1, numerator and denominator of expression (9) all converge to 0. By using L' Hospital rule, we get:

$$F_0^2(z_2)\big|_{z_2=1} = F_0^2(1) = \frac{\mu}{\lambda_1 + \mu} \cdot \lim_{z_2 \to 1} \frac{\lambda_1(1-f(z_2) - \lambda_2(1-z_2)}{\mu/z_2 - [\lambda_1(1-f(z_2) - \lambda_2(1-z_2) + \mu]} p_{000}$$

$$= \frac{\mu}{\lambda_1 + \mu} \cdot \frac{-\lambda_1 f'(z_2) + \lambda_2}{-\frac{\mu}{z_2^2} + \lambda_1 f'(z_2) - \lambda_2}\bigg|_{z_2=1} p_{000} = \frac{\mu \lambda_2}{(\lambda_1 + \mu)(\mu - \lambda_1 - \lambda_2)} p_{000} \quad (12)$$

Derivating both sides of expression (10) again and substituting $f(1), f'(1)$ into expression (10), we have $f''(1) = 2\mu \lambda_2^2 / (\mu - \lambda_1)^3$. After substituting $f(1), f'(1), f''(1)$ into expression (9) and using Mathematica software, we get the following expression:

$$\frac{\partial F_0^2(z_2)}{\partial z_2}\bigg|_{z_2=1} = F_0^{2'}(1) = \frac{\mu \lambda_2 [\lambda_1^3 + \lambda_1^2(3\lambda_2 - \mu) + \mu^3 + \lambda_1(2\lambda_2^2 \mu - \lambda_2 \mu - \mu^2)]}{(\mu + \lambda_1)^2 (\mu - \lambda_1)(\mu - \lambda_1 - \lambda_2)} \quad (13)$$

## 4. Performance evaluation
### 4.1 State probability of server

$F^1(z_1, z_2) = \sum_{i=1}^{\infty} F_i^1(z_2) z_1^i$ , $F^2(z_1, z_2) = \sum_{i=0}^{\infty} F_i^2(z_2) z_1^i, i \geq 0$ are the generating function of the two classes of customers respectively. When $z_1 = 1$ and $z_2 = 1$, $F^1(1,1)$ and $F^2(1,1)$ represent class-1 and class-2 customers' probablity of occupying server respectively. Substituting expression (12) into expression (6) and (7), by using L' Hospital rule again, we get:

$$F^1(1,1) = \left(\frac{\mu - \lambda_2}{\mu - \lambda_1 - \lambda_2} - 1\right) p_{000} \quad , \quad F^2(1,1) = \frac{\lambda_2}{\mu - \lambda_1 - \lambda_2} p_{000}$$

Due to $F^1(1,1) + F^2(1,1) + p_{000} = 1$, we have

$$p_{000} = \frac{\mu - \lambda_1 - \lambda_2}{\mu}.$$

Then we give the following theorem:

**Theorem 1:** In non-preemptive M/M/1 queueing system, clas-1 customer has non-preemptive priority over class-2 customer, then the state probablity of the sever is occupied by calss-1 customer, class-2 customer and the server is free is given by following expressions respectively:

$$P_{class-1} = \frac{\lambda_1}{\mu} \quad , P_{class-2} = \frac{\lambda_2}{\mu} \quad , P_{free} = \frac{\mu - \lambda_1 - \lambda_2}{\mu}$$

### 4.2 average length of two classes of customers
Substituting expression (12) and (13) into expression (6) and (7), by using mathematics software, after complicated deducing, we have :

$$\frac{\partial F^1(z_1,1)}{\partial z_1}\bigg|_{z_1=1} = \frac{\lambda_1^2\lambda_2 + \mu^2\lambda_1}{\mu^2(\mu-\lambda_1)}, \qquad \frac{\partial F^2(z_1,1)}{\partial z_1}\bigg|_{z_1=1} = \frac{\lambda_1\lambda_2}{\mu^2}$$

$$\frac{\partial F^2(1,z_2)}{\partial z_2}\bigg|_{z_2=1} = [(\frac{\lambda_2}{\mu}-1) + \frac{\mu+\lambda_1}{\lambda_2}(1+\frac{\lambda_2^2}{\mu^2}-\frac{\lambda_2}{\mu})]F_0^2(1) + (1-\frac{\lambda_2}{\mu})(\frac{\mu+\lambda_1}{\lambda_2})F_0^{2'}(1)$$

$$\frac{\partial F^2(1,z_2)}{\partial z_2}\bigg|_{z_2=1} = \frac{\lambda_1(\lambda_1+\lambda_2+\mu)}{\mu^2}F_0^2(1) + \frac{\mu+\lambda_1}{\mu}F_0^{2'}(1)$$

**Theorem 2:** In non-preemptive M/M/1 queueing system, class-1 customer has non-preemptive priority over class-2 customers .

The average length of class-1 customer is given by following expression :

$$L_1 = \frac{\partial(F^1(z_1,1)+F^2(z_1,1))}{\partial z_1}\bigg|_{z_1=1}$$

$$= \frac{\lambda_1^2\lambda_2+\mu^2\lambda_1}{\mu^2(\mu-\lambda_1)} + \frac{\lambda_1\lambda_2}{\mu^2} = \frac{\lambda_1(\mu+\lambda_2)}{\mu(\mu-\lambda_1)}$$

The average length of class-2 customer is given by following expression :

$$L_2 = \frac{\partial(F^1(1,z_2)+F^2(1,z_2))}{\partial z_2}\bigg|_{z_2=1}$$

$$= [(\frac{\lambda_1(\mu+\lambda_1+\lambda_2)}{\mu^2}) + \frac{\lambda_2-\mu}{\mu} + \frac{\mu+\lambda_1}{\lambda_2}(1+\frac{\lambda_2^2}{\mu^2}-\frac{\lambda_2}{\mu})]F_0^2(1) + (\frac{\mu+\lambda_1}{\lambda_2})F_0^{2'}(1)$$

$$= \frac{\lambda_1^2\lambda_2+\mu^3+\lambda_1(2\lambda_2^2+\mu^2)}{\mu^2(\lambda_1+\mu)} + \frac{\mu[\mu^3+\lambda_1^3+\lambda_1^2(3\lambda_2-\mu)+\lambda_1(2\lambda_2^2-\mu\lambda_2-\mu^2)]}{(\mu^2-\lambda_1^2)(\mu-\lambda_1-\lambda_2)^2}$$

Here $F_0^2(1), F_0^{2'}(1)$ is given by expression （12）和（13）respectively .

## 5.Conclusion

Discussion in this article enriches modern queueing theory, has wide use in practical questions. As in computer communication net, in order to offer multi-layer quality of service for different kinds of customer, priority control is necessary. Sometimes, number of priority is not only two. For this circumstance, we can use the poisson stream's property of folding, considering one flow of these multi-class customers and one folding flow, all importing customer may be looked as only two classes of customers. Using the conclusion in the article, we can get corresponding result similar to that in the article. These results are useful for optimizing the performance of transfer unit and raising the efficiency of net.